\documentclass[twoside]{article}
\usepackage{amssymb,amsmath}
\usepackage[a4paper,left=1.5in,right=1.5in,top=1in,bottom=1in,headsep=5mm]{geometry}
\usepackage{mathrsfs}
\usepackage{hyperref}

\newtheorem{thm}{Theorem}[section]

\newtheorem{prop}[thm]{Proposition}

\newtheorem{rem}[thm]{Remark}

\numberwithin{thm}{section}
 \numberwithin{equation}{section} \allowdisplaybreaks
\newcommand{\norm}[1]{\left\|#1\right\|}
\newcommand{\nnorm}[1]{\left|\!\left|\!\left|#1\right|\!\right|\!\right|}
\newcommand{\abs}[1]{\left|#1\right|}

\newcommand{\Real}{\mathbb R}

\newcommand{\no}{\nonumber}
\newcommand{\eps}{\varepsilon}

\newcommand{\T}{\mathscr{T}}

\newcommand{\F}{\mathscr{F}}

\renewcommand{\Re}{\mathrm{Re}}
\newcommand{\ls}{\leqslant}
\newcommand{\gs}{\geqslant}
\newcommand{\dif}{\mathrm{d} }

\newenvironment{pf}[1][Proof]{{\par \textbf{#1.}}\;}{\hfill
$\Box$\par}
\begin{document}
 \pagestyle{myheadings}

\markboth{\hfill \textrm{Hao, Hsiao, and Li}
 \hfill}{\hfill \textrm{Well-posedness for Rotating GPE in 3D} \hfill}


\title{\parbox{\textwidth}{\Large\textbf{Global well-posedness for the Gross-Pitaevskii equation
with an angular momentum rotational term in three dimensions}}}%

\author{\begin{minipage}{.87\textwidth}\large Chengchun Hao\thanks{Corresponding author, e-mail:
hcc@amss.ac.cn.}
 and Ling Hsiao\\[1mm]
\normalsize\textsl{Institute of Mathematics, Academy of Mathematics
\& Systems Science, CAS, Beijing 100080, People's Republic of
China}\\[2mm]
\large Hai-Liang Li\\[1mm]
\normalsize\textsl{Department of Mathematics, Capital Normal
University, Beijing 100037, People's Republic of
China}\\[2mm]
(Received \hfill)\end{minipage}}
\date{}
\maketitle

\begin{abstract}
In this paper, we establish the global well-posedness of the Cauchy
problem for the Gross-Pitaevskii equation with an angular momentum
rotational term in which the angular velocity is equal to the isotropic trapping frequency in the space $\Real^3$.

\textbf{Key words:} {Gross-Pitaevskii equation; angular momentum
rotation; harmonic trap potential; global well-posedness.}

\textbf{2000 Mathematics Subject Classification:} {35Q55, 35A05}
\end{abstract}

\section{Introduction}

The Gross-Pitaevskii equation (GPE), derived independently by Gross
\cite{Gro61} and Pitaevskii \cite{Pit61}, arises in various models
of nonlinear physical phenomena. This is a Schr\"odinger-type
equation with an external field potential $V_{ext}(t,x)$ and a local
cubic nonlinearity:
\begin{align}\label{cgpe}
    i\hbar\partial_t u+\frac{\hbar^2}{2m}\Delta u=V_{ext}u+\beta
    \abs{u}^2u.
\end{align}

The GPE \eqref{cgpe} in physical dimensions ($2$ and $3$ dimensions)
is used in the meanfield quantum theory of Bose-Einstein condensate
(BEC) formed by ultracold bosonic coherent atomic ensembles. A rigorous derivation of the GPE for the dynamics of (non-rotating) BEC has been obtained by Erd\"os, Schlein and Yau \cite{ESY06, ESY07}.
Recently, several research groups \cite{HMEWC,MCWD00,MCBD01,MAHHWC}
have produced quantized vortices in trapped BECs, and a typical
method they used is to impose a laser beam on the magnetic trap to
create a harmonic anisotropic rotating trapping potential. Seiringer has discussed the stationary GPE in \cite{Sei02}, Lieb and Seiringer have rigorously derived for the description of the ground state asymptotics of rotating Bose gases in \cite{LS06}. The
properties of BEC in a rotational frame at temperature $T$ being
much smaller than the critical condensation temperature $T_c$
\cite{LL77} are well described by the macroscopic wave function
$u(t,x)$, whose evolution is governed by a self-consistent, mean
field nonlinear Schr\"odinger equation (NLS) in a rotational frame,
also known as the Gross-Pitaevskii equation with an angular momentum
rotation term:
\begin{align}\label{rgpe}
    i\hbar\partial_t u+\frac{\hbar^2}{2m}\Delta u=V(x)u+NU_0
    \abs{u}^2u-\Omega L_z u, \; x\in\Real^3,\, t\gs 0,
\end{align}
where the wave function $u(t,x)$ corresponds to a condensate state,
$m$ is the atomic mass, $\hbar$ is the Planck constant, $N$ is the
number of atoms in the condensate, $\Omega$ is the angular velocity
of the rotating laser beam, and $V(x)$ is an external trapping
potential. When a harmonic trap potential is concerned,
$V(x)=\frac{m}{2}\sum_{j=1}^3\omega_{j}^2x_j^2$ with $\omega_{1}$,
$\omega_{2}$ and $\omega_{3}$ being the trap frequencies in the
$x_1$-, $x_2$- and $x_3$-direction, respectively. The local
nonlinearity term $NU_0\abs{u}^2u$ arises from an assumption about
the delta-shape interatomic potential, where $U_0=4\pi \hbar^2
a_s/m$ describes the interaction between atoms in the condensate
with $a_s$ (positive for repulsive interaction and negative for
attractive interaction) the $s$-wave scattering length.
$L_z=-i\hbar(x_1\partial_{x_2}-x_2\partial_{x_1})$ is the third
component of the angular momentum $L=x\times P$ with the momentum
operator $P=-i\hbar\nabla$.

After normalization, proper nondimensionalization and dimension
reduction in certain limiting trapping frequency regime \cite{PS03},
the system \eqref{rgpe} becomes to the dimensionless GPE in
$d$-dimensions ($d=2,3$):
\begin{align}\label{dlgpe}
    iu_t+\frac{1}{2}\Delta u=V_d(x)u+\beta_d\abs{u}^2u-\Omega L_z
    u,\;x\in\Real^d,\, t>0,
\end{align}
where $L_z=i(x_1\partial_{x_2}-x_2\partial_{x_1})$ and
\begin{align}\label{V}
    \beta_d=\left\{\begin{array}{l}\beta\sqrt{\gamma_3/2\pi},\\
    \beta,\end{array}\right.\quad V_d(x)=\left\{\begin{array}{ll}
    (\gamma_1^2x_1^2+\gamma_2^2x_2^2)/2, \;&d=2,\\
    (\gamma_1^2x_1^2+\gamma_2^2x_2^2+\gamma_3^2x_3^2)/2,
    &d=3,\end{array}\right.
\end{align}
with $\gamma_1>0$, $\gamma_2>0$ and $\gamma_3>0$ constants,
$\beta=\frac{4\pi a_s N}{a_0}$, $a_0=\sqrt{\frac{\hbar}{m\omega_m}}$
and $\omega_m=\min\{\omega_1,\omega_2,\omega_3\}$.

In general, it is a rather complicated process about the dynamics of
solutions (in particular, vortex) for GPE \eqref{rgpe} under the
interaction of trapping frequencies and angular rotating motion. The
recent numerical simulation of GPE \eqref{rgpe} for different choice
of trap frequencies $(\gamma_1,\gamma_2)$ can help us to understand
the complicated dynamical phenomena caused by the angular rotating
and spatial high frequency motion. The case of different frequency
$\gamma_1\neq\gamma_2$ gives much complicated behavior and thus is
rather difficult to be studied rigorously \cite{BDZ06,BWM05}. To our
knowledge, the equation \eqref{rgpe} has been only investigated for
some specific cases by numerical simulation. Therefore, to develop
methods for constructing analytical solutions to the GPE
\eqref{cgpe} or some specific cases is the first step in order to
understand the dynamics caused by the trapping and rotation.

To begin with, we first consider the case
$\gamma_1=\gamma_2=\gamma_3=\omega$ which means the spatial isotropic
motion. In order to derive the exact analytic formula for the solution to the linear equation, we have to assume that the angular velocity is equal to the isotropic trapping frequency, i.e., $\Omega=\omega$. In the present paper, we focus on the Cauchy problem of the Gross-Pitaevskii
equation with an angular momentum rotational term in three
dimensions
\begin{align}\label{gpe}
    &iu_t+\frac{1}{2}\Delta u=\frac{\omega^2}{2}\abs{x}^2 u+\beta\abs{u}^{2}u-\omega L_z
    u,\;
    x\in\Real^3,\, t\gs 0, \\
    &u(0,x)=u_0(x),\, x\in\Real^3,\label{data}
\end{align}
where the wave function $u=u(t,x):
[0,\infty)\times\Real^3\to\mathbb{C}$ corresponds to a condensate
state, $\Delta$ is the Laplace operator on $\Real^3$, $\omega\gs 1$
and $\beta>0$ are constants, and
$L_z=-i(x_1\partial_{x_2}-x_2\partial_{x_1})=i(x_2\partial_{x_1}-x_1\partial_{x_2})$
is the dimensionless angular momentum rotational term.

We assume that the initial value
\begin{align}\label{data.space}
    u_0(x)\in \Sigma:=\{u\in H^1(\Real^3):\, \abs{x}u\in
    L^2(\Real^3)\},
\end{align}
with the norm
\begin{align*}
    \norm{u}_\Sigma =\norm{u}_{H^1}+\norm{\abs{x}u}_{L^2}.
\end{align*}

Note that for multi-dimensional GPE~\eqref{rgpe},
nothing is known about the exact integration except for the case
 $\gamma_1=\gamma_2$ ($=\gamma_3$ for 3D) considered in
\cite{Car02,Car03} without the angular momentum rotational term,
namely, $\Omega=0$. 

In the case of two dimensions \cite{HHL}, the linear operator
$i\partial_t+\frac{1}{2}(\Delta-\omega^2\abs{x}^2+2\omega L_z)$ can
be written as the form $i\partial_t+\frac{1}{2}(\nabla -ib)^2$ where
$b=(-\omega x_2, \omega x_1)$ satisfies the Coulomb gauge condition
$\nabla\cdot b=0$. However, in  three dimensions, we can not find
such a potential $b$ satisfying the Coulomb gauge condition
$\nabla\cdot b=0$. So that it is impossible to write the linear
operator as a similar operator as in the case of two dimensions, and
this makes the problem be much more difficult than the two
dimensional case.

In addition, there are three ingredients that play important roles
in the proof of our result. The first involves the solution of the
Cauchy problem to the linear equation
\begin{align}\label{lin.eqn}
\begin{split}
    &iu_t+\frac{1}{2}\Delta u=\frac{\omega^2}{2}\abs{x}^2 u-\omega L_z
    u,\;
    x\in\Real^3,\, t\gs 0, \\
    &u(0,x)=u_0(x),\, x\in\Real^3,
\end{split}
\end{align}
which is significant for investigating the properties of the
evolution operator corresponding to the linear operator
$i\partial_t+\frac{1}{2}\Delta-\frac{\omega^2}{2}\abs{x}^2+\omega
L_z$. The second one is to obtain the Strichartz estimates for the
foregoing linear operator. The last one is that there exist two
Galilean operators $J(t)$ and $H(t)$ (as blow) which can commute
approximatively with the linear operator and can be viewed as the
substitute of $\nabla$ and $x$ respectively in the non-potential
case.

Now we state the main result of this paper.
\begin{thm}\label{thm}
Let $u_0\in\Sigma$ and $\rho\in(2,6)$. Then, there exists a unique
solution $u(t,x)$ to the Cauchy problem \eqref{gpe}--\eqref{data}.
And the solution satisfies, for any $T\in(0,\infty)$, that
\begin{align*}
    u(t,x), J(t)u(t,x), H(t)u(t,x)\in \mathcal{C}(\Real;
    L^2(\Real^3))\cap L^{\gamma(\rho)}(0,T; L^\rho(\Real^3)),
\end{align*}
where $\frac{2}{\gamma(\rho)}=\frac{3}{2}-\frac{3}{\rho}$, $J(t)$
and $H(t)$ are defined as below as in \eqref{J} and \eqref{H},
respectively.
\end{thm}

\begin{rem}
Since the GPE \eqref{gpe} (or \eqref{dlgpe}) in a rotational frame
is time reversible and time transverse invariant, the above result
is also valid for the case when $t<0$.
\end{rem}

The paper is organized as follows.  In Sec.~\ref{Sec.Str}, the
evolution operator of the linear equation and the Strichartz
estimates about the former operator are first established.
Sec.~\ref{Sec.law} is devoted to the derivation of some conservation
identities such as the mass, the energy, the angular momentum
expectation, and the pseudo-conformal conservation laws in the whole
space $\Real^3$ for \eqref{gpe}--\eqref{data}. Finally, the
nonlinear estimates and the proof of Theorem~\ref{thm} are obtained
in Sec.~\ref{Sec.pf}.

\section{The Strichartz estimates and some main
operators}\label{Sec.Str}

On the analogy of the Schr\"odinger operators with magnetic fields
in \cite{AHS78,Caz03}, we can also define the propagator associated with the self-adjoint operator
\begin{align*}
    H_0\equiv \frac{1}{2}\left[-\Delta+\omega^2\abs{x}^2-2\omega L_z\right],
\end{align*}
i.e., $S(t)=e^{-iH_0t}$, which can be explicitly expressed as
\begin{align}\label{def.S}
    S(t)\phi:=&\left(\frac{\omega}{2\pi i\sin (\omega
    t)}\right)^{\frac{3}{2}}\int_{\Real^3}e^{i\omega\left(\frac{\abs{x-y}^2}{2}\cot(\omega t)-\tilde{x}\cdot
    y\right)} \phi(y)dy\\
    \equiv&\left(\frac{\omega}{2\pi i\sin (\omega
    t)}\right)^{\frac{3}{2}}e^{i\omega \frac{\abs{x}^2}{2}\cot(\omega t)}
    \int_{\Real^3}e^{i\omega A(t)x\cdot y}e^{i\omega \frac{\abs{y}^2}{2}\cot(\omega
    t)}\phi(y)dy,\label{def2.S}\\
    \equiv&\left(\frac{\omega}{2\pi i\sin (\omega
    t)}\right)^{\frac{3}{2}}e^{i\omega \frac{\abs{x}^2}{2}\cot(\omega t)}
    \int_{\Real^3}e^{i\omega x\cdot A^T(t)y}e^{i\omega \frac{\abs{y}^2}{2}\cot(\omega
    t)}\phi(y)dy,\label{def3.S}
\end{align}
for $0<t\ls
\frac{\pi}{4\omega}$ through a complicated computation, where $\tilde{x} :=(-x_2,x_1,(\csc(\omega t)-\cot (\omega t))x_3)$,
the matrix $A(t)$ is defined by
\begin{align*}
             A(t)=\left(
                         \begin{array}{ccc}
                           \cot(\omega t) & -1 & 0 \\
                           1 & \cot(\omega t) & 0 \\
                           0 & 0 & \csc(\omega t) \\
                         \end{array}
                       \right),
\end{align*}
and $A^T(t)$ is the transpose of the matrix $A(t)$. Note
that this formula is valid only for small time, due to the singularity formation for the fundamental
solution.

Next, we derive the dual operator $S^*(t)$ of $S(t)$.  From
\begin{align*}
    &\langle S(t)f,g\rangle=\langle
    \widehat{S(t)f(x)},\overline{\widehat{\bar{g}(x)}}\rangle
    =\langle
    \widehat{S(t)f}(\xi),\widehat{g}(-\xi)\rangle\\
    =&\left\langle\left(\frac{\sin (\omega
    t)}{i\omega}\right)^{\frac{3}{2}}\F\left(e^{i\omega \frac{\abs{x}^2}{2}\cot(\omega
    t)}\right)(\xi)\right.\\
    &\left.\qquad\qquad*\left(
    e^{i \frac{\abs{(A^T(t))^{-1}\xi}^2}{2\omega}\cot(\omega
    t)}f((A^T(t))^{-1}\xi/\omega)\right),\,\widehat{g}(-\xi)\right\rangle\\
    =&\left(\frac{\sin (\omega
    t)}{i\omega}\right)^{\frac{3}{2}}\int_{\Real^3}\int_{\Real^3}
    \F\left(e^{i\omega \frac{\abs{x}^2}{2}\cot(\omega
    t)}\right)(\eta)\\
    &\qquad\qquad\cdot\left(
    e^{i \frac{\abs{(A^T(t))^{-1}(\xi-\eta)}^2}{2\omega}\cot(\omega
    t)}f((A^T(t))^{-1}(\xi-\eta)/\omega)\right)d\eta\widehat{g}(-\xi)d\xi\\
    =&\left(\frac{\sin (\omega
    t)}{i\omega}\right)^{\frac{3}{2}}\int_{\Real^3}\int_{\Real^3}
    \F\left(e^{i\omega \frac{\abs{x}^2}{2}\cot(\omega
    t)}\right)(\eta)\\
    &\qquad\qquad\cdot\left(
    e^{i \frac{\abs{(A^T(t))^{-1}\zeta}^2}{2\omega}\cot(\omega
    t)}f((A^T(t))^{-1}\zeta/\omega)\right)\widehat{g}(-\zeta-\eta)d\zeta
    d\eta\\
    =&\left(\frac{\sin (\omega
    t)}{i\omega}\right)^{\frac{3}{2}}\int_{\Real^3}
    e^{i \frac{\abs{(A^T(t))^{-1}\zeta}^2}{2\omega}\cot(\omega
    t)}f((A^T(t))^{-1}\zeta/\omega)\\
    &\qquad\qquad\qquad\qquad\qquad \cdot\F\left(e^{i\omega \frac{\abs{\cdot}^2}{2}\cot(\omega
    t)}g(-\cdot)\right)(\zeta)d\zeta\\
    =&\left(\frac{\omega}{i\sin (\omega
    t)}\right)^{\frac{3}{2}}\int_{\Real^3}f(\xi) e^{i\omega \frac{\abs{\xi}^2}{2}\cot(\omega
    t)}\F\left(e^{i\omega \frac{\abs{\cdot}^2}{2}\cot(\omega
    t)}g(-\cdot)\right)(\omega A^T(t)\xi)d\xi\\
    =&\langle f, S^*(t)g\rangle,
\end{align*}
we can define the dual operator as
\begin{align}
    S^*(t)g:=&\left(\frac{\omega}{2\pi i\sin (\omega
    t)}\right)^{\frac{3}{2}}e^{i\omega \frac{\abs{x}^2}{2}\cot(\omega
    t)}\int_{\Real^3}e^{-i\omega A^T(t)x\cdot y}e^{i\omega \frac{\abs{y}^2}{2}\cot(\omega
    t)}g(-y)dy\nonumber\\
    =&\left(\frac{\omega}{2\pi i\sin (\omega
    t)}\right)^{\frac{3}{2}}e^{i\omega \frac{\abs{x}^2}{2}\cot(\omega
    t)}\int_{\Real^3}e^{i\omega A^T(t)x\cdot y}e^{i\omega \frac{\abs{y}^2}{2}\cot(\omega
    t)}g(y)dy.\label{dualS}
\end{align}

In order to obtain the Strichartz estimates, we have to estimate the
norm $\norm{S(t)S^*(s)}_{L^1\to L^\infty}$. Indeed, we have the
following proposition:

\begin{prop}
Let $u(t)=S(t)u_0$, then $u(t)$ satisfies the linear equation
\eqref{lin.eqn}. And the operator $S(t)$ has following properties:\\
\indent (1) $S(t)$ is unitary on $L^2$, i.e. $\norm{S(t)}_{L^2\to L^2}=1$;\\
\indent (2) $\norm{S(t)S^*(s)}_{L^1\to L^\infty}\lesssim
\abs{t-s}^{-3/2}$ for $0<s<t\ls \frac{\pi}{4\omega}$.
\end{prop}

\begin{pf} By computation, it is easy to see that
\begin{align}\label{sc.1}
    iu_t=&\left(\frac{\omega}{2\pi i\sin (\omega
    t)}\right)^{\frac{3}{2}}\int_{\Real^3}e^{i\omega\left(\frac{\abs{x-y}^2}{2}\cot(\omega t)-\tilde{x}\cdot
    y\right)}\Big[-\frac{3}{2}i\omega \cot(\omega t)\nonumber\\
    &+\omega^2\csc^2(\omega t)\frac{\abs{x-y}^2}{2}+\omega^2(\csc(\omega t)-\cot(\omega t))\csc(\omega t) x_3y_3)
    \Big] u_0(y)dy,
\end{align}
\begin{align}\label{sc.2}
    \frac{1}{2}\Delta u=&\left(\frac{\omega}{2\pi i\sin (\omega
    t)}\right)^{\frac{3}{2}}\int_{\Real^3}e^{i\omega\left(\frac{\abs{x-y}^2}{2}\cot(\omega t)-\tilde{x}\cdot
    y\right)}\Big[\frac{3}{2}i\omega \cot(\omega t)\nonumber\\
    &-\omega^2\cot^2(\omega t)\frac{\abs{x-y}^2}{2}-\omega^2\frac{\abs{y}^2}{2}+\omega^2\cot(\omega t)\, \tilde{x}
    \cdot y\Big] u_0(y)dy,
\end{align}
\begin{align}\label{sc.3}
    -\frac{\omega^2}{2}\abs{x}^2u=\left(\frac{\omega}{2\pi i\sin (\omega
    t)}\right)^{\frac{3}{2}}\int_{\Real^3}e^{i\omega\left(\frac{\abs{x-y}^2}{2}\cot(\omega t)-\tilde{x}\cdot
    y\right)}\left(-\frac{\omega^2}{2}\abs{x}^2\right)u_0(y)dy,
\end{align}
and
\begin{align}\label{sc.4}
    \omega L_z u=\left(\frac{\omega}{2\pi i\sin (\omega
    t)}\right)^{\frac{3}{2}}&\int_{\Real^3}e^{i\omega\left(\frac{\abs{x-y}^2}{2}\cot(\omega t)-\tilde{x}\cdot
    y\right)}\omega^2\Big[(x_2y_1-x_1y_2)\cot(\omega t)\nonumber\\
    &\quad \quad\quad \quad\quad \quad\quad \quad +x_1y_1+x_2y_2\Big] u_0(y)dy.
\end{align}
Summing \eqref{sc.1}--\eqref{sc.4}, we have
\begin{align*}
    iu_t+\frac{1}{2}\Delta u-\frac{\omega^2}{2}\abs{x}^2 u+\omega L_z
    u=0,
\end{align*}
which yields the desired result.

For the unitariness of $S(t)$ on $L^2$, it is trivial by using the Plancherel
theorem. Thus, we omit the detailed proof.

Next, we prove the dispersive property of $S(t)$. Indeed, we have for
$0<s<t\ls \frac{\pi}{4\omega}$
\begin{align*}
    &S(t)S^*(s)f\\
    =&\left(\frac{\omega^2}{-(2\pi)^2 \sin (\omega
    s)\sin (\omega t)}\right)^{\frac{3}{2}}e^{i\omega\frac{\abs{x}^2}{2}\cot(\omega t)}
    \int_{\Real^3}e^{i\omega A(t)x\cdot y}e^{i\omega\frac{\abs{y}^2}{2}\cot(\omega
    t)}\\
    &\qquad\qquad \cdot e^{i\omega \frac{\abs{y}^2}{2}\cot(\omega s)}\int_{\Real^3}
    e^{i\omega A^T(s)y\cdot z}e^{i\omega \frac{\abs{z}^2}{2}\cot(\omega s)}
    f(z)dzdy\\
    =&\left(\frac{\omega^2}{-2\pi \sin (\omega
    s)\sin (\omega t)}\right)^{\frac{3}{2}}e^{i\omega\frac{\abs{x}^2}{2}\cot(\omega t)}
    \int_{\Real^3}e^{i\omega A(t)x\cdot y}e^{i\omega\frac{\abs{y}^2}{2}[\cot(\omega
    t)+\cot(\omega s)]}\\
    &\qquad\qquad \cdot
    \F^{-1}\left(e^{i\omega \frac{\abs{\cdot}^2}{2}\cot(\omega s)}
    f(\cdot)\right)(\omega A^T(s)y)dy\\
    =&\left(\frac{\sin (\omega
    s)}{-2\pi \sin (\omega t)}\right)^{\frac{3}{2}}e^{i\omega\frac{\abs{x}^2}{2}\cot(\omega t)}
    \int_{\Real^3}e^{i A(t)x\cdot (A^T(s))^{-1}\xi}\\
    &\qquad\qquad \cdot e^{i\frac{\abs{(A^T(s))^{-1}\xi}^2}{2\omega}[\cot(\omega
    t)+\cot(\omega s)]}
    \F^{-1}\left(e^{i\omega \frac{\abs{\cdot}^2}{2}\cot(\omega s)}
    f(\cdot)\right)(\xi)d\xi\\
    =&\left(\frac{\sin (\omega
    s)}{-2\pi \sin (\omega t)}\right)^{\frac{3}{2}}e^{i\omega\frac{\abs{x}^2}{2}\cot(\omega t)}
    \int_{\Real^3}e^{i B(t,s)x\cdot \xi}e^{i\frac{\abs{(A^T(s))^{-1}\xi}^2}{2\omega}[\cot(\omega
    t)+\cot(\omega s)]}\\
    &\qquad\qquad \cdot
    \F^{-1}\left(e^{i\omega \frac{\abs{x}^2}{2}\cot(\omega s)}
    f(x)\right)(\xi)d\xi\\
    =&\left(\frac{\sin (\omega
    s)}{-\sin (\omega t)}\right)^{\frac{3}{2}}e^{i\omega\frac{\abs{x}^2}{2}\cot(\omega t)}
    \Big[\Big(\F^{-1}e^{i\frac{\abs{(A^T(s))^{-1}\xi}^2}{2\omega}[\cot(\omega
    t)+\cot(\omega s)]}\Big)\\
    & \qquad\qquad\qquad\qquad *\big(e^{i\omega \frac{\abs{x}^2}{2}\cot(\omega s)}
    f(-x)\big)\Big](B(t,s)x)\\
    =&\left(\frac{\sin (\omega
    s)}{-\sin (\omega t)}\right)^{\frac{3}{2}}e^{i\omega\frac{\abs{x}^2}{2}\cot(\omega t)}
    \\
    & \cdot
    \int_{\Real^3}\Big(\F^{-1}e^{i\frac{\abs{(A^T(s))^{-1}\xi}^2}{2\omega}[\cot(\omega
    t)+\cot(\omega s)]}\Big)(B(t,s)x-y) e^{i\omega \frac{\abs{y}^2}{2}\cot(\omega s)}
    f(-y)dy
\end{align*}
where the matrix $B(t,s)$ is given by
\begin{align*}
    B(t,x)=\left(
      \begin{array}{ccc}
         B_{11}(t,s)& B_{12}(t,s) & 0 \\
        -B_{12}(t,s) & B_{11}(t,s) & 0 \\
        0 & 0 & \csc(\omega t) \sin(\omega s) \\
      \end{array}
    \right),
\end{align*}
with
\begin{align*}
       B_{11}(t,s)=\sin^2(\omega s)(\cot(\omega t)\cot(\omega s)+1),\;
      B_{12}(t,s)=\sin^2(\omega s)(\cot(\omega t)-\cot(\omega s)).
\end{align*}

 Noticing that
\begin{align*}
    &\abs{\Big(\F^{-1}e^{i\frac{\abs{(A^T(s))^{-1}\xi}^2}{2\omega}[\cot(\omega
    t)+\cot(\omega s)]}\Big)(B(t,s)x-y)}\\
    =&(2\pi)^{-\frac{3}{2}}\abs{\int_{\Real^3}e^{i(B(t,s)x-y)\cdot\xi}e^{i\frac{\abs{(A^T(s))^{-1}\xi}^2}{2\omega}[\cot(\omega
    t)+\cot(\omega s)]}d\xi}\\
    =&\csc^3(\omega
    s)\left(\frac{\omega}{\pi[\cot(\omega t)+\cot(\omega s)]}\right)^{\frac{3}{2}}\abs{\int_{\Real^3}e^{\frac{i\sqrt{2\omega}}{\sqrt{\cot(\omega
    t)+\cot(\omega s)}}(B(t,s)x-y)\cdot A^T(s)\eta}
    e^{i\abs{\eta}^2}d\eta}\\
    \lesssim &\left(\frac{\omega}{\pi\sin^2(\omega
    s)[\cot(\omega t)+\cot(\omega s)]}\right)^{\frac{3}{2}},
\end{align*}
we can get, for $0<s<t\ls \frac{\pi}{4\omega}$, that
\begin{align*}
    \norm{S(t)S^*(s)f}_{L^\infty}\lesssim &\left(\frac{\omega}{\pi\sin(\omega t)\sin(\omega
    s)[\cot(\omega t)+\cot(\omega
    s)]}\right)^{\frac{3}{2}}\int_{\Real^3}\abs{f(y)}dy\\
    = &\left(\frac{\omega}{\pi\sin(\omega
    (t+s))}\right)^{\frac{3}{2}}\norm{f}_{L^1}
    \lesssim  \abs{t+s}^{-3/2}\norm{f}_{L^1}\\
    \lesssim &\abs{t-s}^{-3/2}\norm{f}_{L^1},
\end{align*}
since $\abs{\sin t}\gs \frac{2}{\pi}\abs{t}$ for
$\abs{t}\ls\frac{\pi}{2}$.
\end{pf}

 Thus, we can obtain similar Strichartz
estimates to the linear Schr\"odinger operator
$e^{i\frac{t}{2}\Delta}$ by the standard methods (c.f. \cite{KT98})
provided that only finite time intervals are involved
(c.f.\cite{Car03}).

\begin{prop}\label{pr.Str}
Let $I$ be an interval contained in $[0,\pi/4\omega]$. Then, it
holds that

\mbox{\rm (1)} For any admissible pair $(\gamma(p),p)$ (that is,
$2/\gamma(p)=3(1/2-1/p)$ for $2\ls p<6$), there exists $C_p$ such
that for any $\phi\in L^2$
\begin{align}\label{Str.1}
    \norm{S(t)\phi}_{L^{\gamma(p)}(I;L^p)}\ls C_p\norm{\phi}_{L^2}.
\end{align}

\mbox{\rm (2)} For any admissible pairs $(\gamma(p_1),p_1)$ and
$(\gamma(p_2),p_2)$, there exists $C_{p_1,p_2}$ such that
\begin{align}\label{Str.2}
    \norm{\int_{I\cap\{s<t\}}S(t)S^*(s)F(s)\dif
    s}_{L^{\gamma(p_1)}(I;L^{p_1})}\ls
    C_{p_1,p_2}\norm{F}_{L^{\gamma(p_2)'}(I;L^{p'_2})}.
\end{align}
The above constants are independent of $I\subset[0,\pi/4\omega]$.
\end{prop}

The integral equation reads
\begin{align}\label{e.int}
    u(t)=S(t)u_0-i\beta \int_0^t S(t)S^*(s)\abs{u}^{2}
    u(s) \dif s.
\end{align}

Since the initial data belong to $\Sigma$, we naturally need the
estimates of $\nabla S(t) \phi$ and $xS(t)\phi$. In fact, from
\eqref{def.S}, we can compute and obtain that
\begin{align*}
    \nabla S(t)\phi=i\omega x\cot (\omega t) S(t)\phi-i\omega
    S(t)(x\cot(\omega t)-\breve{x})\phi,
\end{align*}
where $\breve{x}=(-x_2,x_1,(\cot\omega t-\csc\omega t)x_3)$, and
\begin{align*}
    &S(t)\nabla \phi=i\omega (x\cot(\omega t)+\tilde{x})
    S(t)\phi-i\omega \cot (\omega t) S(t)(x \phi),
\end{align*}
which yield
\begin{align}
    \nabla S(t)\phi=&\cos(\omega t) S(t)(\cos(\omega t)\nabla- \sin(\omega
    t)\breve{\nabla})\phi\no\\
    &-i\omega \sin (\omega t) S(t)\left[(\cos(\omega t)x-\sin(\omega t)\breve{x}
    )\phi\right],\label{e1}\\
    xS(t)\phi=&\cos(\omega t) S(t)\left[(\cos(\omega t)x-\sin(\omega t)\breve{x}
    )\phi\right]\no\\
    &-\frac{i}{\omega}\sin (\omega t) S(t)(\cos(\omega t)\nabla- \sin(\omega
    t)\breve{\nabla})\phi,\label{e2}
\end{align}
where $\breve{\nabla}=(-\partial_{x_2},\partial_{x_1},(\cot\omega
t-\csc\omega t)\partial_{x_3})$.

Thus, we have
\begin{align*}
    &S(t)(-i\nabla)\phi\no\\
    =&\left[\omega \sin (\omega t)(\cos(\omega t)x+\sin(\omega t)\breve{x})
    -i\cos(\omega t)(\cos(\omega t)\nabla+ \sin(\omega
    t)\breve{\nabla})\right]S(t)\phi,
\end{align*}
and
\begin{align*}
    &S(t)\omega x\phi\no\\
    =&\left[\omega \cos (\omega t)(\cos(\omega t)x+\sin(\omega t)\breve{x})
    +i\sin(\omega t)(\cos(\omega t)\nabla+ \sin(\omega
    t)\breve{\nabla})\right]S(t)\phi.
\end{align*}

 For convenience, we denote
\begin{align}\label{J}
    J(t)=\omega \sin (\omega t)(\cos(\omega t)x+\sin(\omega t)\breve{x})
    -i\cos(\omega t)(\cos(\omega t)\nabla+ \sin(\omega
    t)\breve{\nabla}),
\end{align}
and the corresponding ``orthogonal'' operator
\begin{align}\label{H}
    H(t)=\omega \cos (\omega t)(\cos(\omega t)x+\sin(\omega t)\breve{x})
    +i\sin(\omega t)(\cos(\omega t)\nabla+ \sin(\omega
    t)\breve{\nabla}),
\end{align}
which will appear in the pseudo-conformal conservation law and play
a crucial role in the nonlinear estimates.

By computation, we can obtain the following commutation relation
\begin{align}\label{commu.rel}
\begin{split}
    &\left[J(t),i\partial_t+\frac{1}{2}\Delta -\frac{\omega^2}{2}\abs{x}^2 +\omega
    L_z\right]=O_J(t),\\
    &\left[H(t),i\partial_t+\frac{1}{2}\Delta -\frac{\omega^2}{2}\abs{x}^2 +\omega
    L_z\right]
    =O_H(t),
\end{split}
\end{align}
where
\begin{align*}
    O_J(t)=\left(0,0,2i\omega \sin(\omega t)(
    \omega \sin(\omega t)x_3-i\cos(\omega t)\partial_{x_3})\right),
\end{align*}
and
\begin{align*}
    O_H(t)=\left(0,0,2i\omega \sin(\omega t)(
    \omega \cos(\omega t)x_3+i\sin(\omega t)\partial_{x_3})\right).
\end{align*}

It is clear that for $0<t\ls\frac{\pi}{4\omega}$
\begin{align}
    \abs{O_J(t)u}=&\abs{\left(0,0,\frac{2i\omega \sin(\omega t)}{2\cos(\omega t)-1}J_3(t)u\right)}
    = \abs{\frac{2\omega \sin(\omega t)}{2\cos(\omega
    t)-1}J_3(t)u}\no\\
    \ls &\frac{2\omega^2 t}{\sqrt{2}-1}\abs{J_3(t)u}
    \lesssim \omega^2 t\abs{J(t)u},\label{OJ}
\end{align}
and
\begin{align}
    \abs{O_H(t)u}=&\abs{\left(0,0,\frac{2i\omega \sin(\omega t)}{2\cos(\omega t)-1}H_3(t)u\right)}
        \lesssim \omega^2 t\abs{H(t)u},\label{OH}
\end{align}
 where $J_3(t)$ and $H_3(t)$ are the third component of the operators $J(t)$ and $H(t)$, respectively.

In addition, denote $M(t)=e^{-i\omega \frac{\abs{x}^2}{2}\tan(\omega
t)}$ and $Q(t)=e^{i\omega \frac{\abs{x}^2}{2}\cot(\omega t)}$, then
\begin{align}\label{JH}
\begin{split}
    J(t)=&-i\cos(\omega t) M(t)(\cos(\omega t) \nabla+\sin (\omega t)
    \breve{\nabla})M(-t),\\
    H(t)=&i\sin(\omega t) Q(t)(\cos(\omega t) \nabla+\sin (\omega t)
    \breve{\nabla})Q(-t).
\end{split}
\end{align}

\section{The conserved quantities}\label{Sec.law}

\begin{prop}\label{pr.cons}
Let $u$ be a solution of the equation \eqref{gpe} with the initial
data $\phi\in\Sigma(\Real^3)$. Then, we have the following conserved
quantities for all $t\gs 0$:

 \mbox{\rm (1)}
The $L^2$-norm:
\begin{align}\label{e.mass}
    \norm{u(t)}_{2}=\norm{u_0}_{2}.
\end{align}

 \mbox{\rm (2)}
The energy for the non-rotating part:
\begin{align}\label{e.non}
    E_0(u)=\frac{1}{2}\norm{\nabla
    u}_2^2+\frac{\omega^2}{2}\norm{x
    u}_2^2+\frac{\beta}{2}\norm{u}_{4}^{4}
    =E_0(u_0).
\end{align}

\mbox{\rm (3)} The angular momentum expectation:
\begin{align}\label{e.mom}
    \langle L_z\rangle(t)=\int_{\Real^2}\bar{u}L_z u \dif x =\langle
    L_z\rangle(0).
\end{align}

\mbox{\rm (4)} The pseudo-conformal conservation law:
\begin{align}\label{e.pcch}
\begin{split}
   &\norm{J(t) u}_2^2+\norm{H(t) u}_2^2 +4\cos(\omega t)(1-\cos(\omega t))
   \left[\norm{\omega x_3u}_2^2+\norm{\partial_{x_3}u}_2^2\right]\\
    &\qquad\qquad +\beta\norm{u}_{4}^{4}
       =2E_0(u_0).
\end{split}
\end{align}
\end{prop}

\begin{pf}
For convenience, we introduce
\begin{align*}
    eq(u):=iu_t+\frac{1}{2}\Delta u-\frac{\omega^2}{2}\abs{x}^2 u-\beta\abs{u}^{2}u+\omega L_z
    u.
\end{align*}

It is clear that \eqref{e.mass} holds by applying the $L^2$-inner
product between $eq(u)$ and $\bar{u}$, and then taking the imaginary
part of the resulting equation.

Since we can use the identity \eqref{e.mom} in the proof of
\eqref{e.non}, we derive \eqref{e.mom} first.  Differentiating
$\langle L_z\rangle(t)$ with respect to $t$, and integrating by
parts, we have
\begin{align*}
    \frac{d\langle L_z\rangle(t)}{dt}
    =&i\int_{\Real^3}\left[\bar{u}_t(x_2\partial_{x_1}u-x_1\partial_{x_2}u)
    +\bar{u}(x_2\partial_{x_1}u_t-x_1\partial_{x_2}u_t)\right] \dif
    x\\
    =&\int_{\Real^3}\left[\overline{-iu_t}(x_2\partial_{x_1}u-x_1\partial_{x_2}u)
    -iu_t(x_2\partial_{x_1}\bar{u}-x_1\partial_{x_2}\bar{u})\right] \dif
    x\\
    =&\int_{\Real^3}\left[\left(\frac{1}{2}\Delta \bar{u}-\frac{\omega^2}{2}\abs{x}^2 \bar{u}
    -\beta\abs{u}^{2}\bar{u}+\omega L_z
    \bar{u}\right)(x_2\partial_{x_1}u-x_1\partial_{x_2}u)\right.\\
    &\qquad\left.+\left(\frac{1}{2}\Delta u-\frac{\omega^2}{2}\abs{x}^2 u-\beta\abs{u}^{2}u+\omega L_z
    u\right)(x_2\partial_{x_1}\bar{u}-x_1\partial_{x_2}\bar{u})\right] \dif
    x\\
    =&\int_{\Real^3}\Re(\Delta
    u(x_2\partial_{x_1}\bar{u}-x_1\partial_{x_2}\bar{u}))
    -\frac{\omega^2}{2}(\abs{x}^2(x_2\partial_{x_1}\abs{u}^2-x_1\partial_{x_2}\abs{u}^2))\\
    &\qquad
    -\beta\Re(\abs{u}^{2}(x_2\partial_{x_1}\abs{u}^2-x_1\partial_{x_2}\abs{u}^2))\dif
    x\\
    =&\frac{\omega^2}{2}\int_{\Real^3}(2x_1x_2\abs{u}^2-2x_2x_1\abs{u}^2))\dif
    x\\
    =&0,
\end{align*}
which yields the desired identity \eqref{e.mom}.

Next, we prove the energy conservation for the non-rotating part
\eqref{e.non}. We consider
\begin{align*}
    \Re (eq(u),u_t)=0,
\end{align*}
where $(\cdot,\cdot)$ denotes the $L^2$-inner product. From the
above, we can get
\begin{align*}
    \int_{\Real^3}\left[\frac{1}{2}\partial_t\abs{\nabla
    u}^2+\frac{\omega^2}{2}\partial_t\abs{xu}^2+\frac{\beta}{2}\partial_t\abs{u}^{4}
    +\frac{\omega}{2}\partial_t(\bar{u}L_zu)\right]\dif x=0,
\end{align*}
which implies the identity \eqref{e.non} with the help of
\eqref{e.mom}.

Finally, the pseudo-conformal conservation law \eqref{e.pcch} can be
easily derived from the definition \eqref{J} and \eqref{H} of the
operators $J(t)$ and $H(t)$ with the help of the energy conservation
for the non-rotating part \eqref{e.non}. We omit the details.
\end{pf}

\section{Nonlinear estimates and the proof of
Theorem~\ref{thm}}\label{Sec.pf}

With the help of \eqref{JH}, we can get
\begin{align*}
    J(t)\abs{u}^{2}u=2\abs{u}^{2}J(t)u-u^2\overline{J(t)u},
\end{align*}
which implies, in view of
$\frac{1}{\rho'}+\eps=\frac{2}{q}+\frac{1}{\rho}$ with
$\max\{0,\frac{2}{\rho}-\frac{2}{3}\}<\eps<\frac{1}{\rho}-\frac{1}{6}$,
$\rho\in (2,6)$ and some $q\in [2,6]$, that
\begin{align*}
    \norm{J(t)\abs{u}^{2}u}_{L^{\left(\frac{\rho}{1-\rho\eps}\right)'}}\ls
    C\norm{u}_{L^q}^{2}\norm{J(t)u}_{L^\rho}.
\end{align*}
From the Sobolev embedding theorem and the H\"older inequality, it
yields
\begin{align*}
    &\norm{J(t)\abs{u}^{2}u}_{L^{\gamma\left(\frac{\rho}{1-\rho\eps}\right)'}(0,T;
    L^{\left(\frac{\rho}{1-\rho\eps}\right)'})}\\
    \ls &
    CT^{1-\frac{3}{2}\eps-\frac{2}{\gamma(\rho)}}\norm{u}_{L^\infty(0,T; H^1)}^{2}
    \norm{J(t)u}_{L^{\gamma(\rho)}(0,T; L^\rho)}.
\end{align*}

Similarly, we have
\begin{align*}
    &\norm{H(t)\abs{u}^{2}u}_{L^{\gamma\left(\frac{\rho}{1-\rho\eps}\right)'}(0,T;
    L^{\left(\frac{\rho}{1-\rho\eps}\right)'})}\\
    \ls &
    CT^{1-\frac{3}{2}\eps-\frac{2}{\gamma(\rho)}}\norm{u}_{L^\infty(0,T; H^1)}^{2}
    \norm{H(t)u}_{L^{\gamma(\rho)}(0,T; L^\rho)},
\end{align*}
and
\begin{align*}
    \norm{\abs{u}^{2}u}_{L^{\gamma\left(\frac{\rho}{1-\rho\eps}\right)'}(0,T;
    L^{\left(\frac{\rho}{1-\rho\eps}\right)'})}\ls
    CT^{1-\frac{3}{2}\eps-\frac{2}{\gamma(\rho)}}\norm{u}_{L^\infty(0,T; H^1)}^{2}
    \norm{u}_{L^{\gamma(\rho)}(0,T; L^\rho)}.
\end{align*}

For convenience, we denote
\begin{align*}
    \nnorm{u}_G:=\norm{u}_G+\norm{J(t)u}_G+\norm{H(t)u}_G,
\end{align*}
where $G$ denotes a normalized space. Thus, we have
\begin{align}\label{non.est}
    &\nnorm{\abs{u}^{2}u}_{L^{\gamma\left(\frac{\rho}{1-\rho\eps}\right)'}(0,T;
    L^{\left(\frac{\rho}{1-\rho\eps}\right)'})}\no\\
    \ls &
    CT^{1-\frac{3}{2}\eps-\frac{2}{\gamma(\rho)}}\norm{u}_{L^\infty(0,T; H^1)}^{2}
    \nnorm{u}_{L^{\gamma(\rho)}(0,T; L^\rho)}.
\end{align}

For any  $\rho\in (2,6)$ and $M\gs 2C\norm{u_0}_\Sigma$, define the
workspace $(\mathcal{D},d)$ as
\begin{align*}
    \mathcal{D}:=\{u:\, \nnorm{u}_{L^\infty(0,T;L^2)\cap L^{\gamma(\rho)}(0,T; L^\rho)}\ls M\},
\end{align*}
with the distance
\begin{align*}
    d(u,v)=\nnorm{u-v}_{L^{\gamma(\rho)}(0,T; L^\rho)}.
\end{align*}
It is clear that $(\mathcal{D},d)$ is a Banach space. Let us
consider the mapping $\T :(\mathcal{D},d)\to (\mathcal{D},d)$
defined by
\begin{align*}
    \T: u(t)\mapsto S(t)u_0-i\beta \int_0^t S(t)S^*(s)\abs{u}^{2}
    u(s) \dif s.
\end{align*}

For $u\in (\mathcal{D},d)$, by the commutation relation
\eqref{commu.rel}, \eqref{OJ},  \eqref{OH}, Proposition~\ref{pr.Str}
and the nonlinear estimate \eqref{non.est}, we obtain
\begin{align}\label{non1}
    \nnorm{\T u}_{L^{\gamma(\rho)}(0,T; L^\rho)}
    \ls &C\norm{u_0}_\Sigma+CT^{1-\frac{3}{2}\eps-\frac{2}{\gamma(\rho)}}\norm{u}_{L^\infty(0,T; H^1)}^{2}
    \nnorm{u}_{L^{\gamma(\rho)}(0,T; L^\rho)}\\
    &\qquad +\norm{O_J(t)u}_{L^{\gamma(\rho)}(0,T; L^\rho)}+\norm{O_H(t)u}_{L^{\gamma(\rho)}(0,T; L^\rho)}\no\\
    \ls &M/2+(CT^{1-\frac{3}{2}\eps-\frac{2}{\gamma(\rho)}}M^{2}+2CT)M\no\\
    \ls &M,
\end{align}
where we have taken $T\in (0,\pi/4\omega]$ so small that
$CT^{1-\frac{3}{2}\eps-\frac{2}{\gamma(\rho)}}M^{2}+2CT\ls 1/2$.
Similar to the above, a straightforward computation shows that it
holds
\begin{align}\label{non2}
    &d(\T u,\T v)\no\\
    \ls &\left[CT^{1-\frac{3}{2}\eps-\frac{2}{\gamma(\rho)}}\left(\norm{u}_{L^\infty(0,T; H^1)}^{2}
    +\norm{v}_{L^\infty(0,T; H^1)}^{2}\right)+2CT\right]
    \nnorm{u-v}_{L^{\gamma(\rho)}(0,T; L^\rho)}\no\\
    \ls &\left[CT^{1-\frac{3}{2}\eps-\frac{2}{\gamma(\rho)}}M^{2}+2CT\right]d(u,v)\no\\
    \ls &\frac{1}{2}d(u,v).
\end{align}
Hence, $\T$ is a contracted mapping from the Banach space
$(\mathcal{D},d)$ to itself. By the Banach contraction mapping
principle, we know that there exists a unique solution
$u\in(\mathcal{D},d)$ to \eqref{gpe}--\eqref{data}. In view of the
conservation laws stated in Proposition \ref{pr.cons}, we can use the standard argument to extend it
uniquely to a solution at the interval $[0,\pi/4\omega]$ which
satisfies, for any $t\in[0,\pi/4\omega]$ and $\rho\in(2,6)$, that
\begin{align*}
    u(t,x), J(t)u(t,x), H(t)u(t,x)\in \mathcal{C}(0,\pi/4\omega;
    L^2(\Real^3))\cap L^{\gamma(\rho)}(0,\pi/4\omega; L^\rho(\Real^3)).
\end{align*}

Then, we can extend the above solution to a global one by
translation. In fact, in order to get the solution in the interval
$(\pi/4\omega, \pi/2\omega]$, we can apply a translation
transformation with respect to the time variable $t$ such that the
initial data $u(\pi/4\omega)$ are replaced by $\tilde{u}(0)$. Let
$\tilde{u}(t,x):=u(t-\pi/4\omega,x)$, then we have from the original
equation with initial data $u(\pi/4\omega,x)$
\begin{align}
    &i\tilde{u}_t+\frac{1}{2}\Delta \tilde{u}=\frac{\omega^2}{2}\abs{x}^2 \tilde{u}+\beta\abs{\tilde{u}}^{2}\tilde{u}
    -\omega L_z \tilde{u},\;
    x\in\Real^3,\, t\gs 0, \label{gpe2}\\
    &\tilde{u}(0,x)=\tilde{u}_0(x):=u(\pi/4\omega,x),\; x\in\Real^3.\label{data2}
\end{align}
In the same way, we can get a solution $\tilde{u}(t,x)$  of
\eqref{gpe2}--\eqref{data2} for $t\in[0,\pi/4\omega]$. It is also a
solution $u(t,x)$ to \eqref{gpe}--\eqref{data} for
$t\in[\pi/4\omega, \pi/2\omega]$ and it is unique. Thus, by an
induction argument with the help of those conserved identities
stated in Proposition~\ref{pr.cons}, we can obtain a global solution
$u(t,x)$ to \eqref{gpe}--\eqref{data} satisfying for any
$T\in(0,\infty)$
\begin{align*}
    u(t,x), J(t)u(t,x), H(t)u(t,x)\in \mathcal{C}(\Real;
    L^2(\Real^3))\cap L^{\gamma(\rho)}(0,T; L^\rho(\Real^3)).
\end{align*}

Therefore, we have completed the proof of the main theorem.

\section*{Acknowledgment}

The authors would like to thank the referees for valuable
comments and suggestions on the original manuscript.
C.C.Hao was partially supported by the Scientific Research Startup
Special Foundation on Excellent PhD Thesis and Presidential Award of
Chinese Academy of Sciences, NSFC (Grant No. 10601061), and the
Innovation Funds of AMSS, CAS of China. L.Hsiao was partially
supported by NSFC (Grant No. 10431060). H.L.Li was partially
supported by NSFC (Grant No. 10431060), Beijing Nova Program and the NCET support of the Ministry of Education of China.

\end{document}